\documentclass[11pt,english]{article}
\usepackage{palatino}
\usepackage[T1]{fontenc}
\usepackage[latin1]{inputenc}
\usepackage{a4}
\usepackage{amssymb}

\makeatletter

\newcommand{\boldsymbol}[1]{\mbox{\boldmath $#1$}}

 \usepackage{amsthm, amsmath, amsfonts, amssymb}
 \theoremstyle{plain}    
 \newtheorem{thm}{Theorem}
 \theoremstyle{plain}    
 \newtheorem{prop}[thm]{Proposition} 

 \theoremstyle{plain}    
 \newtheorem{lem}[thm]{Lemma} 
 \theoremstyle{definition}
 \newtheorem{defn}[thm]{Definition}
 \theoremstyle{remark}
 \newtheorem*{rem*}{Remark}

\usepackage{babel}
\makeatother
\begin{document}

\title{Decomposition of symplectic vector fields with respect to a fibration
in lagrangian tori}

\author{Nicolas Roy%
\footnote{\noindent Address : Geometric Analysis Group, Institut für Mathematik,
Humboldt Universität, Rudower Chaussee 25, Berlin D-12489, Germany.
Email : \texttt{roy@math.hu-berlin.de}%
}~%
\footnote{The author would like to thank Janko Latschev for reading the manuscript.%
}}

\maketitle
\begin{abstract}
Given a fibration of a symplectic manifold by lagrangian tori, we
show that each symplectic vector field splits into two parts : the
first is Hamiltonian and the second is symplectic and preserves the
fibration. We then show an application of this result in the study
of the regular deformations of completely integrable systems.
\end{abstract}

\section*{Introduction}

We would like to begin this introduction by motivating the study of
the fibrations in lagrangian tori of a symplectic manifold $\left(\mathcal{M},\omega\right)$.
Such fibrations naturally arise in the study of completely integrable
systems (CI in short). These are the dynamical systems defined by
a Hamiltonian $H\in C^{\infty}\left(\mathcal{M}\right)$ admitting
a momentum map, i.e. a set $\boldsymbol{A}=\left(A_{1},...,A_{d}\right):\mathcal{M}\rightarrow\mathbb{R}^{d}$
of smooth functions, $d$ being half of the dimension of $\mathcal{M}$,
satisfying $\left\{ A_{j},H\right\} =0$ and $\left\{ A_{j},A_{k}\right\} =0$
for all $j,k:1...d$, and whose differentials $dA_{j}$ are linearly
independent almost everywhere. Then, the Arnol'd-Mineur-Liouville
Theorem \cite{arnold_1,mineur,liouville} insures that in a neighbourhood
of any connected component of any compact regular fiber $\boldsymbol{A}^{-1}\left(a\right)$,
$a\in\mathbb{R}^{d}$, of the momentum map, there exists a fibration
in lagrangian tori along which $H$ is constant. These tori are thus
invariant by the dynamics generated by the associated Hamiltonian
vector field $X_{H}$. 

Despite the {}``local'' character of the Arnol'd-Mineur-Liouville
Theorem, it is tempting to try to glue together these {}``local''
fibrations in the case of \emph{regular} Hamiltonians, i.e. those
for which there exists, near each point of \emph{}$\mathcal{M}$,
a local fibration in invariant lagrangian tori. Unfortunately, this
is not always possible. Some Hamiltonians do not admit any (global)
fibration in lagrangian tori, and some others admit several different
ones (the prototype of degenerate Hamiltonian system is the free particle
moving on the sphere $S^{2}$). Nevertheless, these examples belong
to the non-generic (within the class of regular CI Hamiltonians) class
of \emph{degenerate} Hamiltonians and one can show (see e.g. \cite{roy_5})
that imposing a nondegeneracy condition (like e.g. those introduced
by Arnol'd \cite{arnold_1}, Kolmogorov \cite{kolmogorov}, Bryuno
\cite{bruno} or Rüssmann \cite{russmann}) insures that there exists
a fibration of $\mathcal{M}$ in lagrangian tori along which $H$
is constant, and moreover that it is unique. The genericity of nondegeneracy
conditions motivates the study of fibrations in lagrangian tori $\mathcal{M}\overset{\pi}{\rightarrow}\mathcal{B}$. 

Such a fibration actually gives rise to several natural geometric
structures that we review in the first section. In particular there
exists a natural process of averaging any tensor field in the direction
of the fibers. This process then allows us to prove (Theorem \ref{theo_decomposition_vectors})
that each symplectic vector fields splits into two parts : the first
is hamiltonian and the second is symplectic and preserves the fibration.
As an application of this, we consider in the last section the regular
deformations of CI systems and we show (Theorem \ref{theo_deformation_hamiltonian})
a \emph{Hamiltonian} normal form for these deformations.

First, let us fix some basic notations. We denote by $\mathcal{V}\left(\mathcal{M}\right)$
the space of smooth vector fields on the manifold $\mathcal{M}$.
A symplectic form $\omega$ on $\mathcal{M}$ provides a isomorphism
$\omega:\mathcal{V}\left(\mathcal{M}\right)\rightarrow\Omega^{1}\left(\mathcal{M}\right)$,
also denoted by $\omega$, i.e. $\omega\left(X\right)=\omega\left(X,.\right)$
for each $X\in C^{\infty}\left(\mathcal{M}\right)$. The inverse is
denoted by $\omega^{-1}:\Omega^{1}\left(\mathcal{M}\right)\rightarrow\mathcal{V}\left(\mathcal{M}\right)$.
For each vector field, we denote by $\phi_{X}^{t}$ its flow at time
$t$. Let $O\subset\mathcal{M}$ be any subset. We say that a vector
field $X$ is symplectic (resp. Hamiltonian) in $\mathcal{O}$ if
its associated $1$-form $\omega\left(X\right)$ is closed (resp.
exact) in $\mathcal{O}$. To each Hamiltonian $H\in C^{\infty}\left(\mathcal{M}\right)$
we can associate a vector field $X_{H}=-\omega^{-1}\left(dH\right)$.
Now, given a fibration $\mathcal{M}\overset{\pi}{\rightarrow}\mathcal{B}$,
we say that a vector field $\tilde{X}\in\mathcal{V}\left(M\right)$
is the lift of a vector field $X\in\mathcal{V}\left(\mathcal{B}\right)$
if for each $b\in\mathcal{B}$ and each $m\in\pi^{-1}\left(b\right)$
we have $\pi_{*}\left(\widetilde{X}_{m}\right)=X_{b}$.

\section{Geometrical structures of regular CI systems}

\subsection{The period bundle}

Let $\left(\mathcal{M},\omega\right)$ be a symplectic manifold of
dimension $2d$ and $\mathcal{M}\overset{\pi}{\rightarrow}\mathcal{B}$
a locally trivial fibration in lagrangian tori, whose fibers are denoted
by $\mathcal{M}_{b}=\pi^{-1}\left(b\right)$, $b\in\mathcal{B}$.
The tangent spaces $L_{m}=T_{m}\mathcal{M}_{\pi\left(m\right)}$ of
the fibers form an integrable vector subbundle $L=\bigcup_{m\in\mathcal{M}}L_{m}$
of $T\mathcal{M}$. A theorem due to Weinstein \cite{weinstein1}
insures that each leaf of a lagrangian foliation (not necessarily
a fibration) is naturally endowed with an affine structure. This affine
structure can actually be expressed in a very convenient way (see
e.g. \cite{woodhouse1}) in terms of a linear connection $\nabla$
on the leaf, as follows : 

\begin{prop}
Let $\mathcal{N}$ be a leaf of a lagrangian foliation $L$. Let the
operator $\nabla:\mathcal{V}\left(\mathcal{N}\right)\times\mathcal{V}\left(\mathcal{N}\right)\rightarrow\mathcal{V}\left(\mathcal{N}\right)$
be defined by \[
\nabla_{X}Y=\omega^{-1}\left(\widetilde{X}\lrcorner d\left(\widetilde{Y}\lrcorner\omega\right)\right),\]
where $\widetilde{X}\in\Gamma\left(L\right)$ and $\widetilde{Y}\in\Gamma\left(L\right)$
extend $X$ and $Y$ in $\mathcal{V}\left(\mathcal{M}\right)$ and
are everywhere tangent to $L$. Then $\nabla$ defines a torsion-free
and flat connection on $\mathcal{N}$.
\end{prop}
Accordingly, each fiber $\mathcal{M}_{b}$ is endowed with such a
torsion-free and flat connection. Moreover, since the foliation actually
defines a fibration, the holonomy of $\nabla$ must vanish. Indeed,
for each $b\in\mathcal{B}$, any set of smooth functions $f_{1},...,f_{d}\in C^{\infty}\left(\mathcal{B}\right)$
whose differentials $df_{j}$ are linearly independent near $b$,
provides $d$ Hamiltonian vector fields $X_{f_{1}\circ\pi},...,X_{f_{d}\circ\pi}\in\mathcal{V}\left(\mathcal{M}\right)$
everywhere tangent to the fibers, parallel on each fiber and linearly
independent in a neighbourhood of $\mathcal{M}_{b}$. They thus form
a global parallel frame on $\mathcal{M}_{b}$, implying that the holonomy
of $\nabla$ vanishes and that each fiber $\mathcal{M}_{b}$ is endowed
with a structure of a standard%
\footnote{Here, {}``standard'' means holonomy-free. We recall that on the
torus $\mathbb{T}^{d}$ there exist also exotic affine structures
with non-zero holonomy, such as Nagano-Yagi's one \cite{nagano_yagi}.
Some authors have \cite{curras,curras_molino} shown that such affine
structures can occur on the leaves of certain lagrangian foliations.
Such foliations do of course not define a fibration.%
} affine torus. Actually, if we denote by $\mathcal{V}_{\nabla}\left(\mathcal{M}_{b}\right)$
the $d$-dimensional vector space composed of parallel vector fields,
the argument above shows also that the union $\bigcup_{b\in\mathcal{B}}\mathcal{V}_{\nabla}\left(\mathcal{M}_{b}\right)$
is naturally endowed with a structure of a vector bundle over $\mathcal{B}$.
In the sequel, the following characterization for the space of sections
of the previous bundle will be useful.

\begin{prop}
\label{prop_fibration_lag_1}A vector field $X\in\mathcal{V}\left(\mathcal{M}\right)$
is vertical and parallel on each fiber if and only if its associated
$1$-form $\omega\left(X\right)$ is a pull-back, i.e\[
X\in\Gamma\left(\bigcup_{b\in\mathcal{B}}\mathcal{V}_{\nabla}\left(\mathcal{M}_{b}\right)\right)\Longleftrightarrow\omega\left(X\right)\in\pi^{*}\left(\Omega^{1}\left(\mathcal{B}\right)\right).\]

\end{prop}
\begin{proof}
Let $\alpha=\omega\left(X\right)$ be the associated $1$-form and
$L=\bigcup_{m}L_{m}$ the vertical lagrangian foliation tangent to
the fibers. It is a well-know fact that $\alpha$ is an element of
$\pi^{*}\left(\Omega^{1}\left(\mathcal{B}\right)\right)$ if and only
if both restrictions $\left.\alpha\right|_{L}$ and $\left.d\alpha\right|_{L}$
vanish. Since the foliation is lagrangian, the $1$-form $\left.\alpha\right|_{L}$
vanishes if and only if $X$ is vertical. Moreover, it follows from
the definition of the connection $\nabla$ on the fibers, that $X$
is parallel if and only if $\left.d\alpha\right|_{L}$ vanishes. 
\end{proof}
Now, each fiber $\mathcal{M}_{b}$ is isomorphic to the standard torus
$\mathbb{T}^{d}$. Thus, among the parallel vector fields on $\mathcal{M}_{b}$,
we can consider those whose dynamics is $1$-periodic. We denote this
set by \[
\Lambda_{b}=\left\{ X\in\mathcal{V}_{\nabla}\left(\mathcal{M}_{b}\right)\mid\phi_{X}^{1}=\mathbb{I}\right\} .\]
This discrete subset can easily be shown to be a lattice in $\mathcal{V}_{\nabla}\left(\mathcal{M}_{b}\right)$.
We call it the \emph{period lattice}. We will show that the union
$\Lambda=\bigcup_{b\in\mathcal{B}}\Lambda_{b}$ is a smooth lattice
subbundle of $\bigcup_{b\in\mathcal{B}}\mathcal{V}_{\nabla}\left(\mathcal{M}_{b}\right)$
and we call it the \emph{period bundle}. One way to proceed is to
construct explicitly smooth sections of $\bigcup_{b\in\mathcal{B}}\mathcal{V}_{\nabla}\left(\mathcal{M}_{b}\right)$
which are $1$-periodic. Such sections can be constructed as Hamiltonian
vector fields with Hamiltonian given by the following lemma.

\begin{lem}
\label{lem_actions_1_periodique}Let $\mathcal{O}\subset\mathcal{B}$
be an open set and $\theta$ a symplectic potential in $\tilde{\mathcal{O}}=\pi^{-1}\left(\mathcal{O}\right)$.
Let $b\rightarrow\gamma\left(b\right)$ be a family of cycles depending
smoothly on $b$ and such that $\gamma\left(b\right)\subset\mathcal{M}_{b}$
for all $b\in\mathcal{O}$. Let the function $\xi\in C^{\infty}\left(\mathcal{B}\right)$
be defined by \[
\xi\left(b\right)=\int_{\gamma\left(b\right)}\theta.\]
Then, the vector field $X_{\xi\circ\pi}$ associated to the Hamiltonian
$\xi\circ\pi$ is vertical, parallel and $1$-periodic on each torus
$\mathcal{M}_{b}$ with $b\in\mathcal{O}$. Moreover, for all $b\in\mathcal{O}$,
its trajectories on $\mathcal{M}_{b}$ are homotopic to the cycle
$\gamma\left(b\right)$.
\end{lem}
Such a function $\xi$ is called an \emph{action}. We can then show
the following.

\begin{thm}
\label{theo_period_bundle}The period bundle $\Lambda$ is a smooth
lattice subbundle of $\bigcup_{b\in\mathcal{B}}\mathcal{V}_{\nabla}\left(\mathcal{M}_{b}\right)$.
Moreover, for any contractible subset $\mathcal{O}\subset\mathcal{B}$,
all smooth local sections $X\in\Gamma\left(\mathcal{O},\Lambda\right)$
are Hamiltonian in $\tilde{\mathcal{O}}=\pi^{-1}\left(\mathcal{O}\right)$.
\end{thm}
\begin{proof}
First, the fibration being locally trivial, we can find, in a contractible
subset $\mathcal{O}\subset\mathcal{B}$, $d$ smooth families of cycles
$\gamma_{j}\left(b\right)$, $j=1..d$, forming for each $b\in\mathcal{O}$
a basis of $\pi_{1}\left(\mathcal{M}_{b}\right)$. On the other hand,
a theorem due to Weinstein \cite{weinstein1} implies that there exists
a symplectic potential in a neighbourhood of a fiber $\mathcal{M}_{b}$.
We can always choose a smaller $\mathcal{O}$ such that there exist
a symplectic potential $\theta$ in $\tilde{\mathcal{O}}=\pi^{-1}\left(\mathcal{O}\right)$.
For each $j=1..d$, let $\xi_{j}\in C^{\infty}\left(\mathcal{B}\right)$
be the action function of the previous lemma. Then, the Hamiltonian
vector fields $X_{j}=X_{\xi_{j}\circ\pi}$ are in the lattice $\Lambda_{b}$
for each $b\in\mathcal{O}$. Moreover, these are primitive elements
of the lattice since their trajectories are homotopic to the cycles
$\gamma_{j}\left(b\right)$ which form a basis of $\pi_{1}\left(\mathcal{M}_{b}\right)$.
For the same reason, they are linearly independent. The vector fields
$X_{j}$ thus form a smooth family of basis of $\Lambda_{b}$ for
all $b\in\mathcal{O}$. This shows the first part of the theorem.
The second part follows from the fact that each local section $X\in\Gamma\left(\mathcal{O},\Lambda\right)$
decomposes into $X=\Sigma_{j}c_{j}X_{j}$, where the coefficients
$c_{j}$ are integer constant. This implies that $X$ is Hamiltonian
in $\tilde{\mathcal{O}}$.
\end{proof}
This theorem implies the existence of a natural integer flat connection
on the vector bundle $\bigcup_{b\in\mathcal{B}}\mathcal{V}_{\nabla}\left(\mathcal{M}_{b}\right)$
since the lattices $\Lambda_{b}$ provides a way to relate the spaces
$\mathcal{V}_{\nabla}\left(\mathcal{M}_{b}\right)$ for neighbouring
$b$. This connection, may have non-vanishing holonomy. We call it
\emph{monodromy}, since it coincides obviously with the monodromy
of the fibration in tori (without having regards to the symplectic
structure and to the fact that the tori are lagrangian).%
\footnote{Moreover, the symplectic form $\omega$ provides an isomorphism between
the sections of $\bigcup_{b\in\mathcal{B}}\mathcal{V}_{\nabla}\left(\mathcal{M}_{b}\right)$
and those of $T^{*}\mathcal{B}$. This gives the base space $\mathcal{B}$
a natural structure of an affine space, as was discovered by Duistermaat
\cite{duistermaat}. %
}

\subsection{The torus action bundle}

Our discussion so far shows that given a fibration in lagrangian tori
$\mathcal{M}\overset{\pi}{\rightarrow}\mathcal{B}$, there exists
a natural associated torus bundle acting on it. Indeed, for each $b\in\mathcal{B}$,
the quotient \[
\mathcal{G}_{b}=\mathcal{V}_{\nabla}\left(\mathcal{M}_{b}\right)/\Lambda_{b}\]
 is a Lie group isomorphic to the torus $\mathbb{T}^{d}$. This isomorphism
is not canonical, but it can be realised by choosing a basis of $\Lambda_{b}$.
We will denote the elements of $\mathcal{G}_{b}$ by $\left[X_{b}\right]$,
with $X_{b}\in\mathcal{V}_{\nabla}\left(\mathcal{M}_{b}\right)$,
since they are equivalence classes. Taking the union over all $b$,
we get a torus bundle $\mathcal{G}=\bigcup_{b\in\mathcal{B}}\mathcal{G}_{b}$.
It is a smooth bundle since the period bundle $\Lambda$ is so. We
stress the fact that $\mathcal{G}$ is in general not a principal
bundle since there might not exist any global action of $\mathbb{T}^{d}$
on $\mathcal{G}$, because of the presence of monodromy, which precisely
prevents us from choosing a global basis of $\Lambda$. On the other
hand, there exists a distinguished global section, since each fiber
is a group with a well-defined identity element. 

Although we cannot apply the general theory of connections on principal
bundles, there is a natural way to speak about local parallel sections
of $\mathcal{G}$ over a subset $\mathcal{O}\subset\mathcal{B}$.
These sections are simply local sections $b\rightarrow\left[X_{b}\right]$
of $\mathcal{G}$, with $b\rightarrow X_{b}$ being a smooth local
parallel section of $\bigcup_{b\in\mathcal{B}}\mathcal{V}_{\nabla}\left(\mathcal{M}_{b}\right)$.
We denote the set of local parallel sections by $\Gamma_{\nabla}\left(\mathcal{O},\mathcal{G}\right)$.

\begin{prop}
\label{prop_local_section_torus}For each contractible subset $\mathcal{O}\subset\mathcal{B}$,
the space $\Gamma_{\nabla}\left(\mathcal{O},\mathcal{G}\right)$ is
a Lie group isomorphic to the torus $\mathbb{T}^{d}$.
\end{prop}
\begin{proof}
If $\mathcal{O}$ is contractible, then the monodromy vanishes in
$\mathcal{O}$ and there exist local sections $X_{1},...,X_{d}\in\Gamma\left(\mathcal{O},\Lambda\right)$
with $\left\{ X_{j}\left(b\right)\right\} $ generating the lattice
$\Lambda_{b}$ at each $b\in\mathcal{O}$. To each element $\left(t_{1},...t_{d}\right)\in\mathbb{T}^{d}=\mathbb{R}^{d}/\mathbb{Z}^{d}$,
we associate $\left[X\right]=\left[t_{1}X_{1}+...+t_{d}X_{d}\right]\in\Gamma_{\nabla}\left(\mathcal{O},\mathcal{G}\right)$.
One easily verifies that this provides an isomorphism.
\end{proof}
For each $b$, the group $\mathcal{G}_{b}$ naturally acts on $\mathcal{M}_{b}$
in the following way. \begin{eqnarray*}
\mathcal{G}_{b}\times\mathcal{M}_{b} & \rightarrow & \mathcal{M}_{b}\\
\left(\left[X_{b}\right],m\right) & \rightarrow & \left[X_{b}\right]\left(m\right)=\phi_{X_{b}}^{1}\left(m\right),\end{eqnarray*}
where $X_{b}\in\mathcal{V}_{\nabla}\left(\mathcal{M}_{b}\right)$
is a representative of the class $\left[X_{b}\right]$. One can see
easily that this action is commutative, free, transitive and affine
with respect to Weinstein's connection on $\mathcal{M}_{b}$. Now,
given any section $g\in\Gamma\left(\mathcal{G}\right)$, its restriction
$\left.g\right|_{\mathcal{O}}$ to any contractible subset $\mathcal{O}\subset\mathcal{B}$
is of the form $\left.g\right|_{\mathcal{O}}=\left[X\right]$, where
$X\in\Gamma\left(\mathcal{O},\bigcup_{b\in\mathcal{B}}\mathcal{V}_{\nabla}\left(\mathcal{M}_{b}\right)\right)$.
We can then extend the previous fiberwise action of the $\mathcal{G}_{b}$
to a vertical action of the sections of the toric bundle $\mathcal{G}$
on $\mathcal{M}$ by \begin{eqnarray*}
\Gamma\left(\mathcal{G}\right)\times\mathcal{M} & \rightarrow & \mathcal{M}\\
\left(g,m\right) & \rightarrow & \left[X\right]\left(m\right)=\phi_{X}^{1}\left(m\right),\end{eqnarray*}
where $X\in\Gamma\left(\mathcal{O},\bigcup_{b\in\mathcal{B}}\mathcal{V}_{\nabla}\left(\mathcal{M}_{b}\right)\right)$
for any contractible neighbourhood $\mathcal{O}$ of $b=\pi\left(m\right)$.
This is well-defined since another choice $X^{'}$ of the representative
class of $\left[X\right]$ would differ from $X$ only by an element
of $\Gamma\left(\mathcal{O},\Lambda\right)$ which would provide $\phi_{X^{'}-X}^{1}=\mathbb{I}$.
This action naturally inherits the properties of the fiberwise action,
and an additional property arises when we restrict ourselves to the
parallel sections of $\mathcal{G}$.

\begin{prop}
For any contractible subset $\mathcal{O}\subset\mathcal{B}$, $\Gamma_{\nabla}\left(\mathcal{O},\mathcal{G}\right)$
acts vertically on $\mathcal{M}$ in a symplectic way. 
\end{prop}
We call this action the \emph{toric action} of $\mathcal{G}$ on $\mathcal{M}$.
Even if this action is local, it provides a way to average any tensor
field on \emph{$\mathcal{M}$}. Indeed, according to Proposition \ref{prop_local_section_torus},
$\Gamma_{\nabla}\left(\mathcal{O},\mathcal{G}\right)$ is a compact
Lie group provided $\mathcal{O}\subset\mathcal{B}$ is simply connected.
It is thus endowed with its Haar measure $\mu_{\mathcal{G}}$ and
for any tensor field $T$ of any type on $\mathcal{M}$, we can define
its \emph{vertical average} $\left\langle T\right\rangle $ in the
following way. For each $m\in\mathcal{M}$, we set \[
\left\langle T\right\rangle _{m}=\int_{\Gamma_{\nabla}\left(\mathcal{O},\mathcal{G}\right)}\left(\phi_{X}^{1}\right)_{*}\left(T\right)d\mu_{\mathcal{G}},\]
where $\mathcal{O}\subset\mathcal{B}$ is any contractible neighbourhood
of $b=\pi\left(m\right)$. We can check that the definition does not
depend on the choice of $\mathcal{O}$. Choosing a basis $X_{1},...,X_{d}$
of $\Gamma\left(\mathcal{O},\Lambda\right)$ provides an explicit
expression for the averaged tensor, namely \[
\left\langle T\right\rangle _{m}=\int_{0}^{1}dt_{1}...\int_{0}^{1}dt_{d}\left(\phi_{X_{1}}^{t_{1}}\right)_{*}\circ...\circ\left(\phi_{X_{d}}^{t_{d}}\right)_{*}\left(T\right),\]
which is also independent of the choice of the basis. 

A tensor field $T$ is called \emph{invariant under the toric action
of} $\mathcal{G}$, or simply $\mathcal{G}$\emph{-invariant}, if
for each local parallel section $X\in\Gamma_{\nabla}\left(\mathcal{O},\mathcal{G}\right)$
we have $\left(\phi_{X}^{1}\right)_{*}\left(T\right)=T$, or equivalently
$\mathcal{L}_{X}T=0.$ The following properties can be proved in a
straightforward way.

\begin{prop}
\label{prop_moyenne_action_torique_2}According to these definitions,
we have the following basic properties :
\begin{enumerate}
\item $T$ is $\mathcal{G}$-invariant if and only if $\left\langle T\right\rangle =T$.
\item $\left\langle \left\langle T\right\rangle \right\rangle =\left\langle T\right\rangle $.
\item Each $p$-form $\alpha\in\Omega^{p}\left(\mathcal{M}\right)$ verifies
$\left\langle d\alpha\right\rangle =d\left\langle \alpha\right\rangle $.
\item Let $T$ and $S$ be two tensor fields. If $T$ is $\mathcal{G}$-invariant,
then the contraction $T\lrcorner S$ with respect to any two indices
verifies $\left\langle T\lrcorner S\right\rangle =T\lrcorner\left\langle S\right\rangle $. 
\item In particular, if $X\in\mathcal{V}\left(\mathcal{M}\right)$ is a
vector field and $\alpha=\omega\left(X\right)$ its associated $1$-form,
then we have $\omega\left(\left\langle \alpha\right\rangle \right)=\left\langle X\right\rangle $.
\end{enumerate}
\end{prop}

\section{Decomposition of symplectic vector fields}

The averaging process presented in the previous section provides a
way to decompose any symplectic vector field into the sum of a Hamiltonian
vector field and a symplectic vector field preserving the fibration.
The key step is the following lemma.

\begin{lem}
\label{lem_forme_ferme_moyenne_null_exacte}If $\alpha$ is a closed
$1$-form on $\mathcal{M}$ whose vertical average vanishes, then
it is exact. Moreover, one can choose the primitive $f\in C^{\infty}\left(\mathcal{M}\right)$,
$\alpha=df$, with the property $\left\langle f\right\rangle =0$.
\end{lem}
\begin{proof}
Let us work locally in a contractible subset $\mathcal{O}\subset\mathcal{B}$.
There exists a basis $\left(X_{1},\cdots,X_{d}\right)$ of $\Gamma\left(\mathcal{O},\Lambda\right)$.
Choosing an {}``initial point'' $m\left(b\right)$ depending smoothly
on $b\in\mathcal{O}$, i.e a smooth section of the restricted bundle
$\pi^{-1}\left(\mathcal{O}\right)\overset{\pi}{\rightarrow}\mathcal{O}$,
let us consider the smooth family of cycles $\gamma_{j}\left(b\right)$
consisting of the orbits $t\rightarrow\phi_{X_{j}}^{t}\left(m\left(b\right)\right)$.
The homology classes $\left[\gamma_{j}\left(b\right)\right]$ form
for each $b\in\mathcal{O}$ a basis of $H_{1}\left(\mathcal{M}_{b}\right)$.
On the other hand, since the fibration $\mathcal{M}\overset{\pi}{\rightarrow}\mathcal{B}$
is locally trivial and $\mathcal{O}$ is contractible, the classes
$\left[\gamma_{j}\left(b\right)\right]$ form a basis of the homology
of $\tilde{\mathcal{O}}=\pi^{-1}\left(\mathcal{O}\right)$. 

We then show that for each $j=1..d$ and each $b\in\mathcal{O}$,
one has $\int_{\gamma_{j}\left(b\right)}\left\langle \alpha\right\rangle =\int_{\gamma_{j}\left(b\right)}\alpha$.
Indeed, one has \begin{eqnarray*}
\int_{\gamma_{j}\left(b\right)}\left\langle \alpha\right\rangle  & = & \int_{0}^{1}dt\,\left\langle \alpha\right\rangle \left(X_{j}\right)\circ\phi_{X_{j}}^{t}\left(m\left(b\right)\right)\\
 & = & \int_{0}^{1}dt\, X_{j}\lrcorner\left(\phi_{X_{j}}^{-t}\right)_{*}\left\langle \alpha\right\rangle .\end{eqnarray*}
Moreover, expressing the average $\left\langle \alpha\right\rangle $
in terms of the generators $X_{j}$, one obtains \[
\int_{\gamma_{j}\left(b\right)}\left\langle \alpha\right\rangle =\int_{0}^{1}dt_{1}...\int_{0}^{1}dt_{d}\int_{0}^{1}dt\,\left(\phi_{X_{1}}^{t_{1}}\right)_{*}\circ\cdots\circ\widehat{\left(\phi_{X_{j}}^{t_{j}}\right)}_{*}\circ\cdots\circ\left(\phi_{X_{d}}^{t_{d}}\right)_{*}\left(X_{j}\lrcorner\left(\phi_{X_{j}}^{t_{j}-t}\right)_{*}\alpha\right),\]
where the entry below $\widehat{\,\,}$ has been omitted. Then, we
check with a trivial change of variable that \[
\int_{0}^{1}dt_{j}\int_{0}^{1}dt\left(X_{j}\lrcorner\left(\phi_{X_{j}}^{t_{j}-t}\right)_{*}\alpha\right)=\int_{\gamma_{j}\left(b\right)}\alpha.\]
This implies that $\int_{\gamma_{j}\left(b\right)}\left\langle \alpha\right\rangle =\int_{\gamma_{j}\left(b\right)}\alpha$.

Finally, the hypothesis $\left\langle \alpha\right\rangle =0$ yields
$\int_{\gamma_{j}\left(b\right)}\alpha=0$, where the classes $\left[\gamma_{j}\left(b\right)\right]$
form a basis of the homology of $\tilde{\mathcal{O}}=\pi^{-1}\left(\mathcal{O}\right)$,
as shown before. Since $\alpha$ is closed, this implies that $\alpha$
is exact. Thus, there exists a function $f\in C^{\infty}\left(\tilde{\mathcal{O}}\right)$
such that $\alpha=df$ in $\mathcal{\tilde{O}}$. This function is
unique up to a constant. On the other hand, we deduce from the property
$\left\langle df\right\rangle =d\left\langle f\right\rangle $ and
the hypothesis $\left\langle \alpha\right\rangle =0$ that $\left\langle f\right\rangle $
is a constant function. This allows us to choose the primitive $f$
in an unique way, requiring that $\left\langle f\right\rangle =0$.
This criterion is independent of the choice of the basis $\left(X_{1},...,X_{d}\right)$
and thus allows us to find a primitive $f$ of $\alpha$ globally
defined on $\mathcal{M}$.
\end{proof}
This lemma has the following corollary.

\begin{prop}
\label{prop_vect_sympl_moyenne_nulle_hamiltonien}If $X$ is a symplectic
vector field with vanishing vertical average, i.e $\left\langle X\right\rangle =0$,
then $X$ is Hamiltonian and \[
X=X_{H}\textrm{ with }\left\langle H\right\rangle =0.\]

\end{prop}
\begin{proof}
Indeed, let $\alpha=\omega\left(X,.\right)$ be the closed $1$-form
associated with $X$. According to Proposition \ref{prop_moyenne_action_torique_2},
$\left\langle X\right\rangle =0$ if and only if $\left\langle \alpha\right\rangle =0$.
The previous lemma then implies that $\alpha$ is exact $\alpha=dF$,
with $\left\langle F\right\rangle =0$, i.e $X$ is Hamiltonian $X=X_{H}$,
with $H=-F$.
\end{proof}
Symplectic lifted vector fields and $\mathcal{G}$-invariant ones
are related as shown in the folloing proposition.

\begin{prop}
\label{prop_vecteur_G_inv}If $X\in\mathcal{V}\left(\mathcal{M}\right)$
is a \textbf{$\mathcal{G}$}-invariant vector field, then it is a
lift of a vector field $Y\in\mathcal{V}\left(\mathcal{B}\right)$.

If $\tilde{Y}\in\mathcal{V}\left(\mathcal{M}\right)$ is a symplectic
lift of a vector field $Y\in\mathcal{V}\left(\mathcal{B}\right)$,
then it is $\mathcal{G}$-invariant.
\end{prop}
\begin{proof}
Indeed, for each $b$ let $\mathcal{O}\subset\mathcal{B}$ be a contractible
neighbourhood of $b\in\mathcal{B}$. Since the toric action of \textbf{$\mathcal{G}$}
is transitive, then for each points $m$ and $m'$ belonging to the
fiber $\mathcal{M}_{b}$, there exists a $\left[Z\right]\in\Gamma_{\nabla}\left(\mathcal{G},\mathcal{O}\right)$
such that $m'=\phi_{Z}^{1}\left(m\right)$. On the other hand, the
fact that $X$ is \textbf{$\mathcal{G}$}-invariant implies $X_{m'}=\left(\phi_{Z}^{1}\right)_{*}X_{m}$.
Then, using $\pi\circ\phi_{Z}^{1}=\pi$, one sees that $\pi_{*}X_{m'}=\pi_{*}\left(\phi_{Z}^{1}\right)_{*}X_{m}$
and thus $\pi_{*}X_{m'}=\pi_{*}X_{m}$. This proves that $\left.X\right|_{\mathcal{M}_{b}}$
is a lift of the tangent vector $\pi_{*}X_{m}\in T_{b}\mathcal{B}$.
Finally, it follows from the smoothness of $X$ that it is a lift
of a vector field on $\mathcal{B}$. 

We postpone the proof of the second statement to the next section
where we prove it in the slightly more general case of time-dependent
vector fields (Proposition \ref{prop_vecteur_symplec_et_lift}). 
\end{proof}
We now give the announced theorem of decomposition of symplectic vector
fields.

\begin{thm}
\label{theo_decomposition_vectors}Each symplectic vector field $X\in\mathcal{V}\left(\mathcal{M}\right)$
can be written in an unique way \[
X=X_{1}+X_{2},\]
where 
\begin{itemize}
\item $X_{1}$ is a Hamiltonian vector field, $X_{1}=X_{A}$, with $\left\langle A\right\rangle =0$,
where $\left\langle A\right\rangle $ is the vertical average of the
Hamiltonian $A$.
\item $X_{2}$ is symplectic and is the lift to $\mathcal{M}$ of a vector
field on $\mathcal{B}$.
\end{itemize}
Moreover, $X_{2}$ is simply the vertical average of $X$, i.e. $X_{2}=\left\langle X\right\rangle $.
\end{thm}
\begin{proof}
Let $\alpha=\omega\left(X,.\right)$ be the 1-form associated with
$X$, which is closed since $X$ is symplectic. Let $\alpha_{2}=\left\langle \alpha\right\rangle $
be the vertical average of $\alpha$ and let $\alpha_{1}=\alpha-\alpha_{2}$.
The $1$-forms $\alpha_{1}$ and $\alpha_{2}$ are closed since $d\left\langle \alpha\right\rangle =\left\langle d\alpha\right\rangle $.
Thus, the vector fields $X_{1}$ and $X_{2}$, associated with $\alpha_{1}$
and $\alpha_{2}$, are symplectic. On the other hand, one has $\left\langle \alpha_{1}\right\rangle =0$
and thus $\left\langle X_{1}\right\rangle =0$. According to Proposition
\ref{prop_vect_sympl_moyenne_nulle_hamiltonien}, this implies that
$X_{1}$ is Hamiltonian, $X_{1}=X_{A}$, with $\left\langle A\right\rangle =0$.
Finally, $\left\langle \alpha_{2}\right\rangle =\alpha_{2}$ implies
that $\left\langle X_{2}\right\rangle =X_{2}$. By Proposition \ref{prop_vecteur_G_inv},
it is a lift of a vector field on \emph{}$\mathcal{B}$.

Moreover, the decomposition $X=X_{1}+X_{2}$ is the unique one of
this type. Indeed, suppose that there is a second decomposition $X=X_{1}^{'}+X_{2}^{'}$
with the same properties. Taking the vertical average of both expressions,
we obtain $\left\langle X_{1}+X_{2}\right\rangle =\left\langle X_{1}^{'}+X_{2}^{'}\right\rangle $
and thus $\left\langle X_{2}\right\rangle =\left\langle X_{2}^{'}\right\rangle $.
Now, by Proposition \ref{prop_vecteur_G_inv}, both $X_{2}$ and $X_{2}^{'}$
are $\mathcal{G}$-invariant. It follows that $X_{2}=X_{2}^{'}$ and
thus $X_{1}=X_{1}^{'}$.
\end{proof}
With respect to the fibration, the $X_{2}$ part is {}``trivial''
since its flow preserves the fibration and the $X_{1}$ part is Hamiltonian.
We stress the fact that this result still holds in the presence of
monodromy. This theorem is used in the sequel to show a normal form
theorem for regular deformations of CI systems.

\section{Application : deformations of CI systems}

\subsection{Regular deformations of completely integrable systems}

Let $\left(H_{0},\mathcal{M}\overset{\pi}{\rightarrow}\mathcal{B}\right)$
be a \emph{regular CI system} composed of a fibration in lagrangian
tori $\mathcal{M}\overset{\pi}{\rightarrow}\mathcal{B}$ and a Hamiltonian
$H_{0}\in C^{\infty}\left(\mathcal{M}\right)$ constant along the
fibers. It is well-known since Poincaré's work \cite{poincare_123}
that adding a small perturbation $\varepsilon K$ to the CI Hamiltonian
$H_{0}$ will destroy its CI character and yield chaotic behaviours.
Nevertheless, it is important to investigate the space of all CI Hamiltonians,
since these are the starting point of any perturbation theory, like
the celebrated K.A.M. Theory \cite{kolmogorov,arnold_2,moser_1} which
actually tells us that one can say a lot about the perturbed Hamiltonian
$H_{\varepsilon}=H_{0}+\varepsilon K$ when $\varepsilon$ is small.

A first step towards the understanding of the space of all CI systems,
is to restrict ourselves to regular deformations of regular CI hamiltonians,
i.e smooth families of Hamiltonians $H_{\varepsilon}$ which are CI
and regular for each $\varepsilon$. At this point, we would like
to stress the fact that this does \emph{not} imply that $H_{\varepsilon}$
is constant along the fibers of a family of $\mathcal{M}\overset{\pi_{\varepsilon}}{\rightarrow}\mathcal{B}$
\emph{depending smoothly on $\varepsilon$}. Nevertheless, we conjecture
that is is true for the generic class of \emph{non-degenerate} Hamiltonians.
We refer to \cite{roy_5} for a review of different nondegeneracy
conditions and we will now restrict our study to the following class
of deformations. 

\begin{defn}
\label{def_deformation_reg}Let $\left(H_{0},\mathcal{M}\overset{\pi}{\rightarrow}\mathcal{B}\right)$
be a regular CI system and let $H_{\varepsilon}\in C^{\infty}\left(\mathcal{M}\right)$
be a smooth family of Hamiltonians. We say that $H_{\varepsilon}$
is a \emph{regular deformation} of $H_{0}$ if there exist a smooth
family of functions $I_{\varepsilon}\in\pi^{*}\left(C^{\infty}\left(\mathcal{B}\right)\right)$,
with $I_{0}=H_{0}$, and a smooth family of symplectomorphisms $\phi^{\varepsilon}:\mathcal{M}\rightarrow\mathcal{M}$,
with $\phi^{0}=\mathbb{I}$, such that \[
H_{\varepsilon}=I_{\varepsilon}\circ\phi^{\varepsilon}\]
for all $\varepsilon$.
\end{defn}
For our purposes, we will need to work now with time-dependent vector
fields since each smooth family of diffeomorphisms $\phi^{\varepsilon}$,
with $\phi^{0}=\mathbb{I}$, is the flow at time $\varepsilon$ of
the time-dependent vector field $X_{\varepsilon}$ defined by \[
\frac{d\left(f\circ\phi^{\varepsilon}\left(m\right)\right)}{d\varepsilon}=X_{\varepsilon}\left(f\right)\circ\phi^{\varepsilon}\left(m\right)\]
for each smooth function $f\in C^{\infty}\left(\mathcal{M}\right)$
and each point $m\in\mathcal{M}$. We denote this flow by $\phi_{X_{\varepsilon}}^{\varepsilon}$.
In all the following, all the considered family $\phi^{\varepsilon}$
of diffeomorphisms will implicitly depend smoothly on $\varepsilon$
and satisfy $\phi^{0}=\mathbb{I}$. We refer e.g. to \cite{marsden_ratiu}
for a review of the properties of time-dependent vector fields.

\subsection{Normal form for regular deformations}

The aim of this section is to show Theorem \ref{theo_deformation_hamiltonian}
which insures that, by changing the function $I_{\varepsilon}$, one
may assume that $\phi^{\varepsilon}$ is a Hamiltonian flow. This
result is based on Theorem \ref{theo_decomposition_flow} which states
that any family of symplectomorphisms $\phi^{\varepsilon}$ can be
written as the composition of a Hamiltonian flow with a family of
fiber-preserving symplectomorphisms. Let us first define precisely
these two notions.

\begin{defn}
\label{def_deformation_hamiltonian}A family of symplectomorphisms
$\phi^{\varepsilon}$ is called \emph{Hamiltonian} if its vector field
$X_{\varepsilon}$ is Hamiltonian, $X_{\varepsilon}=X_{A_{\varepsilon}}$,
with $A_{\varepsilon}\in C^{\infty}\left(\mathcal{M}\right)$ depending
smoothly on $\varepsilon$. 
\end{defn}

\begin{defn}
A family of diffeomorphisms $\phi^{\varepsilon}:\mathcal{M}\rightarrow\mathcal{M}$
is called \emph{fiber-preserving} if there exists a family of diffeomorphisms
on the base space $\varphi^{\varepsilon}:\mathcal{B}\rightarrow\mathcal{B}$
such that\[
\pi\circ\phi^{\varepsilon}=\varphi^{\varepsilon}\circ\pi.\]
We say that $\phi^{\varepsilon}$ is \emph{vertical} whenever $\varphi^{\varepsilon}=\mathbb{I}$
for all $\varepsilon$.
\end{defn}
Whenever a vector field on $\mathcal{M}$ is both symplectic and a
lift of a vector field on $\mathcal{B}$, then we have the following
property.

\begin{prop}
\label{prop_vecteur_symplec_et_lift}If $\tilde{Y}_{\varepsilon}\in\mathcal{V}\left(\mathcal{M}\right)$
is symplectic for each $\varepsilon$ and is a lift of a time-dependent
vector field $Y_{\varepsilon}\in\mathcal{V}\left(\mathcal{B}\right)$,
then it is $\mathcal{G}$-invariant and for each tensor field $T$
one has \[
\left\langle \left(\phi_{\tilde{Y}_{\varepsilon}}^{\varepsilon}\right)_{*}T\right\rangle =\left(\phi_{\tilde{Y}_{\varepsilon}}^{\varepsilon}\right)_{*}\left\langle T\right\rangle .\]

\end{prop}
\begin{proof}
Let denote by $\phi^{\varepsilon}=\phi_{\tilde{Y}_{\varepsilon}}^{\varepsilon}$
the flow of $\tilde{Y}_{\varepsilon}$. This flow is fiber-preserving
and thus verifies $\pi\circ\phi^{\varepsilon}=\varphi^{\varepsilon}\circ\pi$
with $\varphi^{\varepsilon}:\mathcal{B}\rightarrow\mathcal{B}$ a
family of diffeomorphisms. One can easily show that $\varphi^{\varepsilon}$
is actually the flow of $Y_{\varepsilon}$. 

First of all, for each vertical and parallel vector field $X\in\Gamma\left(\bigcup_{b\in\mathcal{B}}\mathcal{V}_{\nabla}\left(\mathcal{M}_{b}\right)\right)$,
one has $\phi_{*}^{\varepsilon}X\in\Gamma\left(\bigcup_{b\in\mathcal{B}}\mathcal{V}_{\nabla}\left(\mathcal{M}_{b}\right)\right)$.
Indeed, according to Proposition \ref{prop_fibration_lag_1}, $\phi_{*}^{\varepsilon}X$
is vertical and parallel if and only if the $1$-form $\omega\left(\phi_{*}^{\varepsilon}X\right)$
is a pull-back. Now, one has $\omega\left(\phi_{*}^{\varepsilon}X\right)=\left(\left(\phi^{\varepsilon}\right)^{-1}\right)^{*}\left(\omega\left(X\right)\right)$
since $\phi^{\varepsilon}$ is symplectic for each $\varepsilon$.
On the other hand, $\omega\left(X\right)=\pi^{*}\beta$ with $\beta\in\Omega^{1}\left(\mathcal{B}\right)$,
since by hypothesis $X$ is vertical and parallel. Thus, one has \[
\omega\left(\phi_{*}^{\varepsilon}X\right)=\left(\left(\phi^{\varepsilon}\right)^{-1}\right)^{*}\pi^{*}\beta=\pi^{*}\left(\left(\varphi^{\varepsilon}\right)^{-1}\right)^{*}\beta.\]
This proves that $\omega\left(\phi_{*}^{\varepsilon}X\right)$ is
a pull-back and thus $\phi_{*}^{\varepsilon}X$ is vertical and parallel. 

If in addition $X\in\Gamma\left(\Lambda,\mathcal{O}\right)$, with
$\mathcal{O}\subset\mathcal{B}$ a subset, i.e. $X$ is $1$-periodic
in $\pi^{-1}\left(\mathcal{O}\right)$, then so is $\phi_{*}^{\varepsilon}X$
in $\phi^{\varepsilon}\left(\pi^{-1}\left(\mathcal{O}\right)\right)$.
Now, the smooth bundle $\Lambda$ has discrete fibers and $\phi_{*}^{\varepsilon}X$
depends smoothly on $\varepsilon$. This implies that for all $\varepsilon$,
one has $\phi_{*}^{\varepsilon}X=\phi_{*}^{\varepsilon=0}X$ and thus
$\phi_{*}^{\varepsilon}X=X$. Then, the derivative with respect to
$\varepsilon$ shows that $\left[\tilde{Y},X\right]=0$, i.e. $\tilde{Y}$
is $\mathcal{G}$-invariant. By linearity, this is true as well for
all $X\in\Gamma\left(\bigcup_{b\in\mathcal{B}}\mathcal{V}_{\nabla}\left(\mathcal{M}_{b}\right)\right)$. 

Therefore, for each $X\in\Gamma\left(\bigcup_{b\in\mathcal{B}}\mathcal{V}_{\nabla}\left(\mathcal{M}_{b}\right)\right)$
and each $\varepsilon$, $\phi^{\varepsilon}$ commutes with the flow
$\phi_{X}^{t}$. This implies that $\phi^{\varepsilon}$ commutes
with the toric action of $\mathcal{G}$ and thus with the averaging
process, i.e. \[
\left\langle \left(\phi_{\tilde{Y}_{\varepsilon}}^{\varepsilon}\right)_{*}T\right\rangle =\left(\phi_{\tilde{Y}_{\varepsilon}}^{\varepsilon}\right)_{*}\left\langle T\right\rangle \]
for any tensor field $T$. 
\end{proof}
We can now give the following decomposition theorem for families of
symplectomorphisms.

\begin{thm}
\label{theo_decomposition_flow}Each family of symplectomorphisms
$\phi^{\varepsilon}$ decomposes in a unique way as follows : \[
\phi^{\varepsilon}=\Phi^{\varepsilon}\circ\phi_{Z_{\varepsilon}}^{\varepsilon},\]
 where 
\begin{itemize}
\item $\Phi^{\varepsilon}$ is a fiber-preserving family of  symplectomorphisms.
\item $Z_{\varepsilon}=X_{G_{\varepsilon}}$ is an time-dependent Hamiltonian
vector field with $\left\langle G_{\varepsilon}\right\rangle =0$.
\end{itemize}
Moreover, the vector field of $\Phi^{\varepsilon}$ equals to the
average $\left\langle X_{\varepsilon}\right\rangle $, where $X_{\varepsilon}$
is the vector field of $\phi^{\varepsilon}$. 
\end{thm}
\begin{proof}
Let $X_{\varepsilon}$ be the vector field of $\phi^{\varepsilon}$.
Theorem \ref{theo_decomposition_vectors} insures that for each $\varepsilon$,
$X_{\varepsilon}$ decomposes into $X_{\varepsilon}=\tilde{Y}_{\varepsilon}+W_{\varepsilon}$,
where $\tilde{Y}_{\varepsilon}$ is a lift of a vector field $Y_{\varepsilon}\in\mathcal{V}\left(\mathcal{B}\right)$
and $W_{\varepsilon}$ is Hamiltonian. Moreover, by looking more carefully
at the proof of Theorem \ref{theo_decomposition_vectors}, one can
easily check that $\tilde{Y}_{\varepsilon}$ and $W_{\varepsilon}$
depend smoothly on $\varepsilon$, since $\tilde{Y}_{\varepsilon}$
is nothing but the vertical average of $X_{\varepsilon}$.

Let $\Psi^{\varepsilon}$ be the family of symplectomorphisms defined
by $\phi_{\tilde{Y}_{\varepsilon}+W_{\varepsilon}}^{\varepsilon}=\phi_{\tilde{Y}_{\varepsilon}}^{\varepsilon}\circ\Psi^{\varepsilon}$
and let $Z_{\varepsilon}$ be its vector field. On the one hand, $\Phi^{\varepsilon}=\phi_{\tilde{Y}_{\varepsilon}}^{\varepsilon}$
is fiber-preserving since $\tilde{Y}_{\varepsilon}$ is a lift of
a vector field on $\mathcal{B}$. On the other hand, on can check
in a straightforward way that the vector field $X_{\varepsilon}^{3}$
of a composition of flows $\phi_{X_{\varepsilon}^{1}}^{\varepsilon}\circ\phi_{X_{\varepsilon}^{2}}^{\varepsilon}$
is given by the formula $X_{\varepsilon}^{3}=X_{\varepsilon}^{1}+\left(\phi_{X_{\varepsilon}^{1}}^{\varepsilon}\right)_{*}X_{\varepsilon}^{2}$.
Therefore, in our case we have $\tilde{Y}_{\varepsilon}+W_{\varepsilon}=\tilde{Y}_{\varepsilon}+\phi_{\tilde{Y}_{\varepsilon}}^{\varepsilon}\left(Z_{\varepsilon}\right)$
and thus \[
Z_{\varepsilon}=\left(\phi_{\tilde{Y}_{\varepsilon}}^{\varepsilon}\right)_{*}^{-1}\left(W_{\varepsilon}\right).\]
According to Theorem \ref{theo_decomposition_vectors}, $W_{\varepsilon}$
is Hamiltonian and verifies $\left\langle W_{\varepsilon}\right\rangle =0$.
First, this insures that $Z_{\varepsilon}$ is Hamiltonian. Second,
Proposition \ref{prop_vecteur_symplec_et_lift} implies that \[
\left\langle Z_{\varepsilon}\right\rangle =\left(\phi_{\tilde{Y}_{\varepsilon}}^{\varepsilon}\right)_{*}^{-1}\left\langle W_{\varepsilon}\right\rangle =0\]
since $\tilde{Y}_{\varepsilon}$ is symplectic and a lift of a vector
field on $\mathcal{B}$.

Finally, we show that this decomposition is unique. Indeed, suppose
that we have a second decomposition $\phi_{X_{\varepsilon}}^{\varepsilon}=\phi_{\tilde{Y}_{\varepsilon}^{'}}^{\varepsilon}\circ\phi_{Z_{\varepsilon}^{'}}^{\varepsilon}$
with the same properties. The vector field $\tilde{Y}_{\varepsilon}^{'}\textrm{ }$
must be a lift of a vector field on $\mathcal{B}$ since $\phi_{\tilde{Y}_{\varepsilon}^{'}}^{\varepsilon}$
is fiber-preserving. On the other hand, as we mentionned before, we
have the relation $\tilde{X}_{\varepsilon}=\tilde{Y}_{\varepsilon}^{'}+\phi_{\tilde{Y}_{\varepsilon}^{'}}^{\varepsilon}\left(Z_{\varepsilon}^{'}\right)$.
Arguing as before, we can show that $\phi_{\tilde{Y}_{\varepsilon}^{'}}^{\varepsilon}\left(Z_{\varepsilon}^{'}\right)$
is a Hamiltonian vector field with vanishing vertical average. Now,
Theorem \ref{theo_decomposition_vectors} tells us that the decomposition
$X_{\varepsilon}=\tilde{Y}_{\varepsilon}+W_{\varepsilon}$ is unique
and thus $\tilde{Y}_{\varepsilon}^{'}=\tilde{Y}_{\varepsilon}$ and
$Z_{\varepsilon}^{'}=Z_{\varepsilon}$.
\end{proof}
As an application, the following theorem gives a normal form for regular
deformations of a given regular CI system. 

\begin{thm}
\label{theo_deformation_hamiltonian}Let $\left(H_{0},\mathcal{M}\overset{\pi}{\rightarrow}\mathcal{B}\right)$
a regular CI system. If $H_{\varepsilon}$ is a regular deformation
of $H_{0}$, then there exist a family of functions $I_{\varepsilon}\in\pi^{*}\left(C^{\infty}\left(\mathcal{B}\right)\right)$
and a family of Hamiltonian symplectomorphisms $\phi_{X_{G_{\varepsilon}}}^{\varepsilon}$,
with $\left\langle G_{\varepsilon}\right\rangle =0$ such that\[
H_{\varepsilon}=I_{\varepsilon}\circ\phi_{X_{G_{\varepsilon}}}^{\varepsilon}\]
for each $\varepsilon$.
\end{thm}
\begin{proof}
By definition, $H_{\varepsilon}$ is a regular deformation of $H_{0}$
if there exist a family of functions $J_{\varepsilon}\in\pi^{*}\left(C^{\infty}\left(\mathcal{B}\right)\right)$
and a family of symplectomorphisms $\phi^{\varepsilon}$ such that
$H_{\varepsilon}=J_{\varepsilon}\circ\phi^{\varepsilon}$. On the
other hand, Theorem \ref{theo_decomposition_flow} insures that $\phi^{\varepsilon}$
decomposes into $\phi^{\varepsilon}=\Phi^{\varepsilon}\circ\phi_{X_{G_{\varepsilon}}}^{\varepsilon}$,
where $\Phi^{\varepsilon}$ is fiber-preserving and $\left\langle G_{\varepsilon}\right\rangle =0$.
Therefore, we have $H_{\varepsilon}=I_{\varepsilon}\circ\phi_{X_{G_{\varepsilon}}}^{\varepsilon}$,
where the function $I_{\varepsilon}=J_{\varepsilon}\circ\Phi^{\varepsilon}$
is indeed an element of $\pi^{*}\left(C^{\infty}\left(\mathcal{B}\right)\right)$
since $\Phi^{\varepsilon}$ is fiber-preserving.
\end{proof}
\begin{rem*}
We can show \cite{roy_1,roy_5} that the families $I_{\varepsilon}$
and $\phi_{X_{G_{\varepsilon}}}^{\varepsilon}$ in the previous theorem
are actually unique provided we assume that $H_{0}$ is \emph{non-degenerate}.
Nondegeneracy conditions are those used in K.A.M. theories, like for
example those introduced by Arnol'd \cite{arnold_1}, Kolmogorov \cite{kolmogorov},
Bryuno \cite{bruno} or Rüssmann \cite{russmann}. They are all open
conditions on the Hamiltonians in $\pi^{*}\left(C^{\infty}\left(\mathcal{B}\right)\right)$.
Imposing Rüssmann's condition (the weakest one) on $H_{0}$ implies
that the fibration in lagrangian tori along which $H_{0}$ is constant
is unique. This allows to show the uniqueness of the families $I_{\varepsilon}$
and $\phi_{X_{G_{\varepsilon}}}^{\varepsilon}$.
\end{rem*}
\bibliographystyle{plain}
\bibliography{/home/roy/math/biblio/biblio_nico}

\end{document}